\documentclass[twoside,leqno,twocolumn]{article}  
\usepackage{ltexpprt}

\usepackage{amsmath}
\usepackage{amssymb}
\usepackage{amsfonts}
\usepackage{url}
\usepackage{graphicx}
\usepackage{nicefrac}
\usepackage{units}
\usepackage{url}

\newcommand{\E}{\ensuremath{\mathbb{E}}}
\newcommand{\V}{\ensuremath{\mathrm{Var}}}

\newcommand{\Prob}{\ensuremath{\mathbb{P}}}

\newcommand{\R}{\ensuremath{\mathbb{R}}}
\newcommand{\N}{\ensuremath{\mathbb{N}}}
\newcommand{\bo}{\mathrm{O}}

\newenvironment{comment}{\setbox0=\vbox\bgroup}{\egroup}

\def\parsec{\par\noindent}

\def\med{\medskip\parsec}

\newtheorem{thm}{Theorem}[section]
\newtheorem{prop}[thm]{Proposition}

\newtheorem{lem}[thm]{Lemma}

\def\lam{\lambda}
\def\dlam{\dot{\lam}}
\def\ddlam{\ddot{\lam}}

\def\bmu{\mbox{\boldmath $\mu$}}
\def\vnu{\mbox{\boldmath $\nu$}}
\def\vgamma{\mbox{\boldmath $\gamma$}}

\def\vpi{\mbox{\boldmath $\pi$}}
\def\VPsi{\mbox{\boldmath $\psi$}}
\def\bv{{\bf v}}

\def\la{\langle}
\def\ra{\rangle}

\begin{document}
\title{\bf  Towards More Realistic Probabilistic Models for\\ 
Data Structures: 
The External Path Length  in \\ Tries 
under the Markov Model}
\author{Kevin Leckey and Ralph Neininger\\Institute for Mathematics\\
J.W.~Goethe University Frankfurt\\
60054 Frankfurt am Main\\
Germany\\
\{leckey, neiningr\}{@}math.uni-frankfurt.de
   \and Wojciech Szpankowski
\thanks{This author's contribution was  made while visiting
J.W.~Goethe University Frankfurt a.M. with an Alexander von Humboldt
research award. This author was also supported by the NSF
Science and Technology Center on Science of Information Grant
CCF-0939370, NSF Grants DMS-0800568 and CCF-0830140, AFOSR Grant FA8655-
11-1-3076, NSA Grant H98230-08-1-0092, and MNSW grant N206
369739. He is also Visiting Professor at ETI,  Gdansk University of
Technology, Poland.}\\
Department of Computer Science\\ 
Purdue University\\
W. Lafayette, IN 47907\\
 U.S.A.\\
 spa{@}cs.purdue.edu} 

\maketitle

\begin{abstract}
\small\baselineskip=9pt
Tries are among the most versatile and widely used data structures on words. 
They are pertinent to the (internal) structure of (stored) words 
and several splitting procedures used in diverse contexts
ranging from document taxonomy to IP addresses lookup, from data
compression (i.e., Lempel-Ziv'77 scheme) to dynamic hashing, 
from partial-match queries to speech recognition, from leader 
election algorithms to distributed hashing tables and graph compression.
While the performance of tries  under a realistic probabilistic model is of 
significant importance, its analysis, even for simplest memoryless sources, 
has proved difficult. Rigorous findings about inherently complex  
parameters  were rarely analyzed
(with a few notable exceptions) 
under more realistic models of string generations.
In this paper we meet these challenges:
By a novel use of the contraction method combined with
 analytic techniques we prove a central limit 
theorem for the external path length of a trie under
a general Markov source.
In particular, our results  apply to the Lempel-Ziv'77 code. We 
envision that the methods described here will have further 
applications to other trie parameters and data structures.
\end{abstract}


\section{Introduction}

\begin{comment}
Tries are one of the most popular data structures on words.
They are pertinent to (internal) structure of (stored) words \cite{Gusfield97}
and several splitting procedures \cite{kps96} used in diverse contexts
ranging from document taxonomy to IP addresses lookup, from data
compression (i.e., Lempel-Ziv'77 scheme) to dynamic hashing,
from partial-match queries to speech recognition, from leader
election algorithms to distributed hashing tables
\cite{Gusfield97,Knuth98,Mahmoud92,spa-book})
Therefore, understanding tries behavior under realistic probabilistic 
model is of real importance (e.g., to study popular Lempel-Ziv'77 
data compression scheme).
Furthermore, developing a widely 
applicable and simple(r) technique to analyze tries for such models 
(i.e., beyond memoryless sources) is imperative. 
There are many trie parameters of interest, but 
the path length (sum of paths from the root to all 
external nodes) is one of the most arduous to analyze. 
In this paper, we meet all of these challenges:
we prove a central limit theorem for the {\it path length} in a trie for
{\it Markovian sources}  by a novel use of the {\it contraction method}.
\end{comment}

We study the external path length of a trie built over $n$ binary
strings generated by a Markov source. More precisely, we assume
that the input is a sequence of $n$ independent and identically
distributed random strings, each being composed of an infinite
sequence of symbols such that the next symbol depends on the previous
one and this dependence is governed by a given transition matrix
(i.e., Markov model).

Digital trees, in particular, tries have been intensively studied for the 
last thirty years 
\cite{clflva01,dev82,dev84,dev92,dev02,dev05,jare88,jasz89,jasz95,jaszta01,kipr88,kiprsz89,kiprsz94,
Knuth98,Mahmoud92,spa-book}, mostly under Bernoulli (memoryless)  
model assumption.
The typical depth under Markovian model was analyzed in \cite{jasz89,jaszta01}. Size, external path length and height under  more general {\em dynamical sources} were studied in the seminal paper of
Cl{\'e}ment,  Flajolet, and Vall{\'e}e \cite{clflva01}, where in particular asymptotic expressions for expectations are identified as well as the asymptotic distributional behavior of the height, see also \cite{clfiflva09}. For further analysis of tries for probabilistic models beyond Bernoulli (memoryless) sources see  Devroye \cite{dev84,dev92}. 

With respect to Markovian models, to the best of our knowledge, no asymptotic distributions for  the external path length have been derived so far. It is well
known \cite{spa-book} that the external path length is more challenging due to 
stronger dependency. In fact, this is already observed for tries under
Bernoulli model \cite{spa-book}.
In this paper we establish the central limit theorem for the
external path length in a trie built over a Markov model using a novel
use of the {\it contraction method}.

Let us  first briefly review the contraction method.
It was introduced in 1991 by Uwe R\"osler \cite{Ro91} for the 
distributional analysis of the complexity of the Quicksort algorithm. 
Over the last 20 years this approach, which is based on exploiting an 
underlying contracting map on a space of probability distributions, 
has been developed as a fairly universal tool for the analysis of 
recursive algorithms and data structures. Here, randomness may come
from a stochastic model for the input or from randomization within the 
algorithms itself (randomized algorithms). General  developments of this 
method were presented in 
\cite{Ro92,RaRu95,Ro99,NeRu04,NeRu04b,FiKa04,DrJaNe08,JaNe08,NeSu12} 
with numerous applications in Theoretical Computer Science.
 
The contraction method has been used in the analysis of tries and other digital trees only under the 
symmetric Bernoulli model (unbiased memoryless source) \cite[Section 5.3.2]{NeRu04}, 
where limit laws for the size and the external path length of tries were
re-derived. The application of the method there was heavily based on the fact 
that precise expansions of the expectations were available, in particular 
smoothness properties of periodic functions appearing in the linear terms 
as well as  bounds on error terms which were $\bo(1)$ for the size and 
$\bo(\log n)$ for the path lengths. Let us observe that 
even in the asymmetric Bernoulli model such error terms seem to be out of 
reach for classical analytic methods; see the discussion in 
Flajolet, Roux, and Vall{\'e}e \cite{flrova10}. Hence, 
for the more general Markov source model considered in the present paper 
we develop a novel use of the contraction method. 

Furthermore, the contraction method 
applied to Markov sources hits another snag, namely,  
the Markov model is not preserved when decomposing the trie into 
its left and right subtree of the root. The initial distribution of the 
Markov source is changed when looking at these subtrees. 
To overcome these problems a couple of new ideas are used for 
setting up the contraction method: First of all, we will use a system of 
distributional recursive equations, one for each subtree.
We then apply the contraction method to this system of recurrences capturing the 
 subtree processes and prove
normality for   the path lengths conditioned on the initial distribution.
In fact, our approach avoids dealing with multivariate recurrences and instead
we reduce the whole analysis to a system of  one-dimensional equations. 
A comparison of a multivariate approach and our new version with systems of recurrences is 
drawn in Section \ref{rem_multivariate}.

We also need asymptotic expansions of the mean and the variance for applying the 
contraction method. 
However, in contrast to very precise information on periodicities of 
linear terms for the symmetric Bernoulli model  mentioned above  our convergence proof   
does only require the leading order term together with a Lipschitz continuity 
property for the error term. 

In this extended abstract we develop the use of systems of recursive distributional 
equations in the context of the contraction method for the external 
path length of tries under a general Markov source model. 
In particular, we prove the central limit theorem for the external path length, 
a result that had been wanting since Lempel-Ziv'77 code was devised in 1977.
The methodology used is general enough to 
cover related quantities and structures as well. We are confident that our 
approach also applies with  minor adjustments at least to the size of tries, 
 the path lengths of digital search trees and PATRICIA tries under 
the Markov source model as well as other more complex data structures on 
words such as suffix trees.\\

\noindent
{\bf Notations:}
Throughout this paper we use the Bachmann-Landau symbols, in particular the big $\bo$ notation. We declare $x\log x:=0$ for $x=0$, where $\log x$ denotes the natural logarithm.  By $B(n,p)$ with $n\in\N$ and $p\in [0,1]$ the binomial distribution is denoted, by $B(p)$ the Bernoulli distribution with success probability $p$, by ${\cal N}(0,\sigma^2)$ the centered normal distribution 
with variance $\sigma^2>0$. We use $C$ as a generic constant that may change from one occurrence to another.

\section{Tries and the Markov source model}
{\bf The Markov source:} We assume binary data strings over the alphabet $\Sigma=\{0,1\}$ generated by
a homogeneous Markov chain.  In general, a homogeneous Markov chain is  
given by its initial distribution $\mu=\mu_0 \delta_0 + \mu_1 \delta_1$ on $\Sigma$ 
and the transition matrix $(p_{ij})_{i,j \in \Sigma}$. Here, $\delta_x$ denotes the Dirac
 measure in $x\in \R$. Hence, the initial state is $0$ with probability $\mu_0$ and $1$ 
with probability $\mu_1$. We have $\mu_0,\mu_1\in [0,1]$ and $\mu_0+\mu_1=1$. 
A transition from state $i$ to $j$ happens with probability $p_{ij}$, $i,j\in \Sigma$. 
Now, a data string is generated as the sequence of states visited by the Markov chain.  
In the Markov source model assumed subsequently  all data strings are independent 
and  identically distributed according to the given Markov chain.

We always assume that $p_{ij}>0$ for all $i,j\in \Sigma$. Hence, the Markov chain is ergodic and has a stationary distribution, denoted by $\pi=\pi_0 \delta_0 + \pi_1 \delta_1$. We have
\begin{align} \label{stat_dist}
\pi_0=\frac{p_{10}}{p_{01}+p_{10}},\qquad \pi_1=\frac{p_{01}}{p_{01}+p_{10}}.
\end{align}
Note however, that our Markov source model does not require the Markov chain to start in its stationary distribution.  

The case $p_{ij}=1/2$ for all $i,j\in\Sigma$ is essentially the symmetric Bernoulli model (only the first bit may have  a different (initial) distribution). The symmetric Bernoulli model has already been studied thoroughly also with respect to the external path length of tries, see \cite{jare88a,kiprsz89,NeRu04}. 
  It behaves differently compared to the asymmetric Bernoulli model and the other Markov source models, as the variance of the external path length is linear with a periodic prefactor  in the symmetric Bernoulli model. In our cases we will find a larger variance of the order $n\log n$ in Theorem \ref{thm:variance} below.
  We exclude the symmetric Bernoulli model case subsequently.
For later reference, we summarize our  conditions as:
\begin{equation}\label{cond_prob}
  \begin{aligned}
p_{ij}&\in (0,1) \mbox{ for all } i,j\in\Sigma, \\
p_{ij}&\neq \frac{1}{2} \mbox{ for some } (i,j)\in\Sigma^2.
  \end{aligned}
\end{equation}
The entropy rate of the Markov chain plays an important role in the asymptotic 
behavior of tries. In particular, it determines  leading order constants of  
parameters of tries that are related to depths of leaves and its external path length.  
The entropy rate for our Markov chain is given by
\begin{eqnarray}\label{def_ent}
   H:= -\sum_{i,j\in \Sigma} \pi_i\, p_{ij} \log p_{ij}=\sum_{i\in \Sigma} \pi_i H_i, 
\end{eqnarray}
where $H_i:=-\sum_{j\in\Sigma} p_{ij} \log p_{ij}$ is the entropy of a transition from state $i$ to the next state. Thus, $H$ is obtained as weighted average of  the entropies of all possible transitions  with weights according to the stationary distribution $\pi$. \\

\noindent
{\bf Tries:} For a given set of data strings over the alphabet $\Sigma=\{0,1\}$ with each data string a unique infinite path in the infinite complete 
rooted binary tree is associated by identifying left branches 
with bit $0$ and right branches with bit $1$. Each string is stored in 
the unique node on its infinite path that is closest to the root and does not belong to any other 
data path, cf. Figure 1.
It is the minimal prefix of a string that distinguishes this string from all
others; for details see 
the monographs of Knuth  \cite{Knuth98}, Mahmoud \cite{Mahmoud92} or Szpankowski \cite{spa-book}.

\begin{figure}[h]
\includegraphics[scale=0.6]{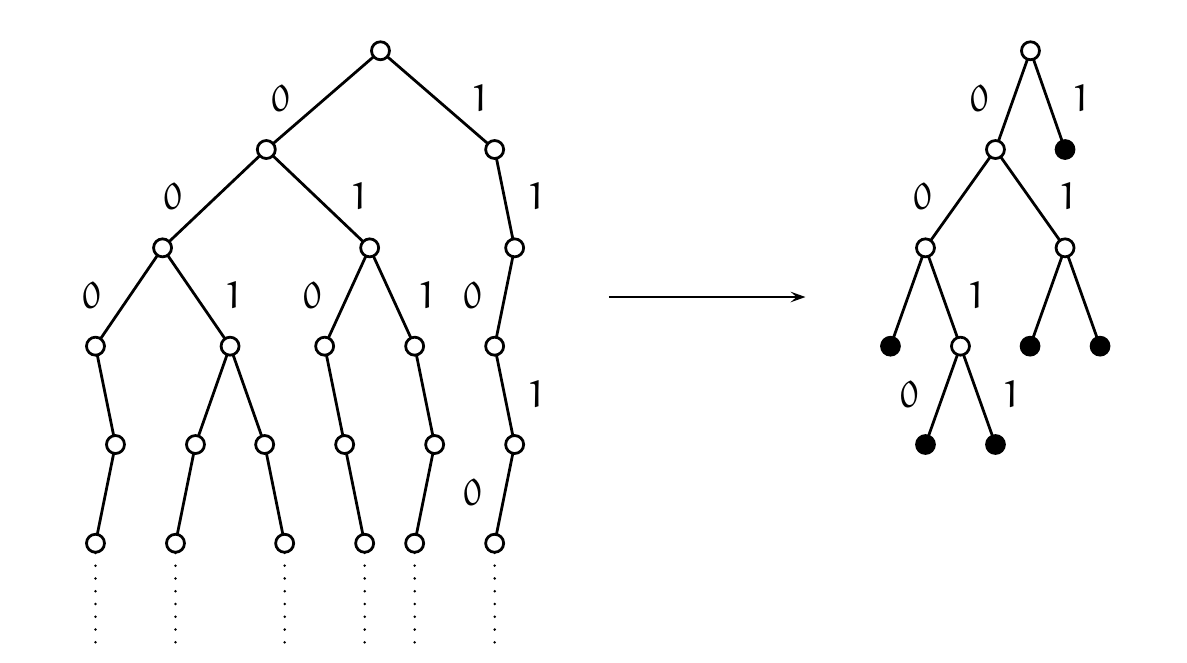}
\caption{The infinite rooted binary tree contains the infinite paths of six strings (left). The corresponding trie
is obtained by cutting each path at the closest node to the root that does not belong to any other path.}
\end{figure}

\section{Recursive Distributional Equations}

For the Markov source model a challenge is to set the right framework under which data
structures to analyze. We formulate in this section a system of distributional recurrences to
capture the distribution of the external path length of tries. Our subsequent analysis is entirely based on these equations. 

We denote by $L_n^\mu$ the external path length of a trie under the Markov source model with initial distribution $\mu$ holding $n$ data. We have $L_0^\mu=L_1^\mu=0$ for all initial distributions $\mu$. The transition matrix is given in advance and suppressed  in the notation. We abbreviate $L_n^i:=L_n^{\delta_i}$ for $i\in\Sigma$. Hence,  $L_n^i$ refers to $n$ independent strings all starting with bit $i$ and then following the Markov chain. We will study $L_n^0$ and $L_n^1$. From the asymptotic behavior of these two sequences we can then directly obtain corresponding results for $L_n^\mu$ for an arbitrary initial distribution $\mu=\mu_0 \delta_0 + \mu_1 \delta_1$ as follows: We denote by $K_n$ the number of data among our $n$ strings which start with bit $0$. Then $K_n$ has the binomial $B(n,\mu_0)$ distribution.
The contributions of the two subtrees of the trie to its external path length 
can be represented by the following stochastic recurrence
\begin{eqnarray} \label{decomp_initial}
L_n^\mu \stackrel{d}{=} L_{K_n}^0 + L_{n-K_n}^1, \qquad n\ge 2,
\end{eqnarray}
where $\stackrel{d}{=}$ denotes that left and right hand side have identical distributions 
and we have that $(L_0^0,\ldots,L_n^0)$, $(L_0^1,\ldots,L_n^1)$ and $K_n$ are independent. We will see later that we can directly transfer asymptotic results for $L_n^0$ and $L_n^1$ to general $L_n^\mu$ via (\ref{decomp_initial}), see, e.g., the proof of Theorem \ref{thm:limit}.

For a recursive decomposition of $L_n^0$ note that we have initial distribution $\delta_0$, thus all data strings start with bit $0$ and are inserted into the left subtree of the root. We denote the root of this left subtree by $w$. 
At node $w$ the data strings are split according to their second bit. We denote by $I_n$ the number of data strings having $0$ as their second bit, i.e., the number of strings being inserted into the left subtree of $w$. The Markov source model implies that $I_n$ is binomial $B(n,p_{00})$ distributed. The right subtree of node $w$ then holds the remaining $n-I_n$ data strings. Consider the left subtree of $w$ together with its root $w$. Conditioned on its number $I_n$ of data strings inserted it is generated by the same Markov source model as the original trie. However, the right subtree of $w$ together with its root $w$ conditioned on its number $n-I_n$ of data strings is generated by a Markov source model with the same transition matrix but another initial distribution, namely $\delta_1$. 
Moreover, by the independence of data strings within the Markov source model, 
these two subtrees are independent conditionally on $I_n$. 
Phrased in a recursive distributional equation we have  
\begin{eqnarray} \label{decomp_zero}
L_n^0 \stackrel{d}{=} L^0_{I_n} + L^1_{n-I_n}+n, \qquad n\ge 2,
\end{eqnarray}
with $(L_0^0,\ldots,L_n^0)$, $(L_0^1,\ldots,L_n^1)$ and $I_n$  independent. A similar arguments yields a recurrence for $L_n^1$. Denoting by $J_n$ a binomial $B(n,p_{11})$ distributed random variable, we have 
\begin{eqnarray} \label{decomp_one}
L_n^1 \stackrel{d}{=} L^0_{n-J_n} + L^1_{J_n}+n, \qquad n\ge 2,
\end{eqnarray}
with $(L_0^0,\ldots,L_n^0)$, $(L_0^1,\ldots,L_n^1)$ and $J_n$ independent.
Our asymptotic analysis of $L_n^\mu$ is based on the distributional 
recurrence system (\ref{decomp_zero})--(\ref{decomp_one})  
as well as (\ref{decomp_initial}).

\section{Analysis of the Mean} \label{sec:mean}

First we study the asymptotic behavior of the expectation of the external path length 
with a precise error term needed to derive a limit law in Section \ref{sec:limit}. The leading order term in Theorem \ref{thm:mean} below has already been derived (even for more general models) in Cl{\'e}ment,  Flajolet and Vall{\'e}e \cite{clflva01}.

\begin{thm}\label{thm:mean}
For the external path length $L_n^\mu$ of a binary trie under the Markov source model with conditions (\ref{cond_prob}) we have 
\begin{align*}
 \E[L_n^\mu] =\frac{1}{H} n\log n + \bo(n), \qquad (n\to\infty),
\end{align*}
with the entropy rate $H$ of the Markov chain given in (\ref{def_ent}). The $\bo(n)$ error term is uniform in the initial distribution $\mu$.
\end{thm}
Our proof of Theorem \ref{thm:mean} as well as the corresponding limit law in Theorem \ref{thm:limit} depend on refined properties of the $\bo(n)$ error term that are first obtained for the initial distributions $\mu=\delta_0$ and $\mu=\delta_1$ and then generalized to arbitrary initial distribution via (\ref{decomp_initial}). For $\mu=\delta_0$ and $\mu=\delta_1$ we
 denote this error term for all $n\in\N_0$ and $i\in \Sigma$ by 
\begin{align}\label{def_fi}
 f_i(n) := \E[L_n^i]-\frac{1}{H} n\log n.
\end{align}
The following Lipschitz continuity of $f_0$ and $f_1$ is crucial for our further analysis:
\begin{prop}\label{thEwLip}
There exists a constant $C>0$ such that for both $i\in\Sigma$ and all $m,n\in\N_0$
\begin{align*}
|f_i(m)-f_i(n)| \leq C |m-n|.
\end{align*}
\end{prop}
The proof of Proposition  \ref{thEwLip} is based on  a refined 
analysis of transfers from growth of toll functions in systems of recursive equations to the growth of the 
quantities itself. The heart of the proof of Proposition  \ref{thEwLip} and hence 
Theorem \ref{thm:mean} is the following transfer result. The proof is technical and provided in the 
full paper version of this extended abstract.
\begin{lem}
\label{lemO1}
Let $(a_i(n))_{n\geq 0}$ and  $(\eta_i(n))_{n\geq 0}$ be real sequences 
and  $(X_{i,n})_{n\geq 2}$  sequences of binomial $B(n,p_i)$ distributed random with $p_i\in(0,1)$ for $i\in\Sigma$.
Assume that for  constants $c_0,c_1,d_0,d_1 \in (0,1)$ with $c_0+d_1=c_1+d_0=1$ we have 
for all  $n\geq 2$ and $i\in\Sigma$ 
\begin{align}
\label{rekfolge}
 a_i(n)&=c_i \E[a_i(X_{i,n})] + d_i \E[a_{1-i}(n-X_{i,n})]\\
& \;\;\;  + \eta_i(n). \nonumber
\end{align}
If furthermore $\eta_i(n)=\bo(n^{-\alpha})$ for an $\alpha>0$ and  both $i\in\Sigma$, then, as $n\to\infty$,
\begin{align*}
 a_i(n)=\bo(1), \quad i\in \Sigma.
\end{align*}
\end{lem}

\section{Analysis of the Variance}

To formulate an asymptotic expansion of the variance of the external path length we  denote by $\lambda(s)$ the largest eigenvalue of the matrix 
$P(s):= (p_{ij}^{-s})_{i,j\in\Sigma}$. Note that $\lambda$ as a function 
of $s$ is smooth. 
We denote its first and second derivative  by $\dlam$ and $\ddlam$ respectively.
Then we have:
\begin{thm} \label{thm:variance}
For the external path length $L_n^\mu$ of a binary trie under the 
Markov source model with conditions (\ref{cond_prob}) we have, as $n\to\infty$,
\begin{align}\label{var_asympt}
\V(L_n^\mu) =\sigma^2 n\log n + o(n\log n), 
\end{align}
where $\sigma^2>0$ is independent of the initial distribution $\mu$ and 
given by 
\begin{equation}
\label{var_formular}
\sigma^2=\frac{\ddlam(-1)-\dlam^2(-1)}{\dlam^3(-1)}.
\end{equation}
With $H_0$ and $H_1$ defined in (\ref{def_ent}) we have 
\begin{align*}
 \sigma^2=&\frac{\pi_0 p_{00}p_{01}}{H^3} \left( \log \left(\frac{p_{00}}{p_{01}}\right)+ \frac {H_1-H_0}{p_{01}+p_{10}}\right)^2\\
&~+\frac{\pi_1 p_{10}p_{11}}{H^3}  \left( \log \left(\frac{p_{10}}{p_{11}}\right)+ \frac {H_1-H_0}{p_{01}+p_{10}}\right)^2.
\end{align*}
\end{thm}

We start with the analysis of the Poisson variance of the external path length, 
i.e. $\tilde{v}_i (\lambda) := \V(L_{N_\lambda}^i)$, $i\in\Sigma$, 
where $N_\lambda$ has the Poisson($\lambda$) distribution and is independent of 
$(L_n^i)_{n\geq 0}$. In the second part we use depoissonization techniques of 
\cite{js98} to obtain the asymptotic behavior of $\V(L_n^i)$. 

The reason why we consider a Poisson number of strings is that for $N_\lambda$ 
i.i.d.~strings with initial distribution $\delta_i$ the number $N_{\lambda p_{i0}}$ of 
strings whose second bit equals 0 and the number $M_{\lambda p_{i1}}$ of 
strings whose second bit equals 1 are independent and remain Poisson distributed. 
Hence, in the Poisson case we obtain similarly to (\ref{decomp_zero}) 
and (\ref{decomp_one}) that for $i\in\Sigma$
\begin{align} \label{decomp_pois}
L_{N_\lambda}^i &\stackrel{d}{=} L_{N_{\lambda p_{i0}}}^0 + L_{M_{\lambda p_{i1}}}^1 \\
& \quad +N_{\lambda p_{i0}}+M_{\lambda p_{i1}}-\textbf{1}_{ \{N_{\lambda p_{i0}}+M_{\lambda p_{i1}}=1\}} \nonumber
\end{align}
where $(L_n^0)_{n\geq 0}$, $(L_n^1)_{n\geq 0}$, $N_{\lambda p_{i0}}$ and 
$M_{\lambda p_{i1}}$ are independent, $N_{\lambda p_{i0}}$ has 
Poisson($\lambda p_{i0}$) distribution and $M_{\lambda p_{i1}}$ has 
Poisson($\lambda p_{i1}$) distribution. Note that 
$\textbf{1}_{ \{N_{\lambda p_{i0}}+M_{\lambda p_{i1}}=1\}}$ is necessary 
in order that (\ref{decomp_pois}) holds when $\{N_\lambda=1\}$.

We denote by $\tilde{\nu}_i (\lambda):= \E [ L_{N_\lambda}^i ]$, 
$i\in\Sigma$, the Poisson expectation of the external path length which is
$$\tilde{\nu}_i (\lambda)=\sum_{n=0}^\infty e^{-\lambda} 
\frac{\lambda^n}{n!} \E [L_n^i].$$
Note that (\ref{decomp_pois}) implies 
\begin{align}
 \label{decomp_pois_exp}
\tilde{\nu}_i(\lambda)=\tilde{\nu}_0 (\lambda p_{i0}) + \tilde{\nu}_1 
(\lambda p_{i1} ) +  \lambda (1- e^{-\lambda}).
\end{align}
We need precise information about the mean (second order term) to derive
the leading term of the variance. We shall use analytic techniques, namely
the Mellin transform as surveyed in \cite{spa-book} that we discuss next.
A Mellin transform
$f^*(s)$ of a real function $f(x)$ is defined as
$$
f^*(s)=\int_0^\infty f(x) x^{s-1} dx.
$$
Let $\nu^*_i(s)$ be the Mellin transform of $\tilde{\nu}_i(\lambda)$. Then,
by known properties of the Mellin transform \cite{spa-book}, the
functional equation (\ref{decomp_pois_exp}) becomes an algebraic  equation
for $i\in \Sigma$
$$
\nu^*_i(s)=\Gamma(s+1)+p_{i0}^{-s}\nu^*_0(s)+p_{i1}^{-s}\nu^*_1(s).
$$
Define the column vector $\vnu^*(s):=(\nu^*_0(s), \nu^*_1(s))$ and the
column vector $\vgamma(s):=(\Gamma(s), \Gamma(s))$.
Then we can write the latter equations as the matrix equation 
$\vnu^*(s)=\vgamma(s+1)+P(s)\vnu^*(s)$ that we write as
\begin{equation}
\label{spa1}
\vnu^*(s)=(I-P(s))^{-1}\vgamma(s+1).
\end{equation}
Then the Mellin transform $\nu^*(s)$ of the mean external path length 
$\E[L^\mu_{N_\lambda}]$ under
the Poisson model satisfies
\begin{equation}
\label{spa2}
\nu^*(s)=\Gamma(s+1)+\bmu(s)\vnu^*(s)
\end{equation}
where $\bmu(s):=(\mu_0^{-s}, \mu_1^{-s})$.

To recover the mean external path length under the Poisson model we 
need to apply the
singularity analysis to (\ref{spa2}).   
For matrix $P(s)$, we define the principal left eigenvector $\vpi(s)$,
the principal right eigenvector $\VPsi(s)$ associated with the
largest eigenvalue $\lambda(s)$ 
such that $\la \vpi(s), \VPsi(s)\ra =1$ where
we write  $\langle {\bf x}, {\bf y}\rangle$ for the
inner product of vectors ${\bf x}$ and ${\bf y}$.
  Then by the {\it spectral representation}
\cite{spa-book} of
$P(s)$  we find
$$
\nu^*(s)=\frac{\Gamma(s)\VPsi(s)}{1-\lambda(s)} + o(1/(1-\lambda(s)))
$$
that leads to the following asymptotic expansion around $s=-1$
\begin{align}
\label{spa3}
\nu^*(s)&= \frac{-1}{\dlam(-1)}\frac{1}{(s+1)^2}+\frac{1}{s+1}\left(\frac{\gamma}{\dlam(-1)}+\frac{\ddlam(-1)}{2\dlam(-1)}\right) \\
& \quad +\frac{1}{s+1}\left(-\frac{\langle \dot\bmu(-1)\dot\VPsi(-1)\rangle}{\dlam(-1)} +1 \right) +\bo(1) \nonumber
\end{align}
where $\dot{\bf x}(t)$ and
$\ddot{\bf x}(t)$ denote the first and second derivatives of the 
vector ${\bf x}(t)$ at $t$. 

Using (\ref{spa3}), inverse Mellin transform, and the residue theorem of Cauchy,
as well as analytic depoissonization of Jacquet and Szpankowski \cite{js98}
we finally obtain
\begin{align}
\label{spa4}
\E[L^\mu_n] &= \frac 1H n \log n + n \left( \frac{\gamma}{\dlam(-1)}+\frac{\ddlam(-1)}{2\dlam(-1)}\right) \\
& \quad +n\left(-\frac{\langle \dot\bmu(-1)\dot\VPsi(-1)\rangle}{\dlam(-1)} +1 +
\Phi(\log n)\right)
+o(n) \nonumber
\end{align}
where $\Phi(x)$ is a periodic function of small amplitude under 
certain rationality condition (and zero otherwise); see \cite{jaszta01} for
details.

The asymptotic analysis of the variance  follows the same pattern, however, it is  more
involved.
Our analysis of the Poisson variance $\tilde{v}_i(\lambda)={\rm Var}(L^i_{N_\lambda})$ 
is based on the following decomposition:

\begin{lem}\label{lem_decomp_pois_var} For any $\lambda>0$ and $i\in\Sigma$ we have
\begin{align}
 \label{decomp_pois_var}
\tilde{v}_i (\lambda) &= \tilde{v}_0 (\lambda p_{i0}) + \tilde{v}_1 
(\lambda p_{i1}) + 2\lambda p_{i0} \tilde{\nu}_0^\prime 
(\lambda p_{i0}) \\
&\quad + 2 \lambda p_{i1} \tilde{\nu}_1^\prime (\lambda p_{i1})
+ 2 \lambda e^{-\lambda} (\tilde{\nu}_0 (\lambda p_{i0})+\tilde{\nu}_1 (\lambda p_{i1}))\nonumber \\
 & \quad + \lambda(1-e^{-\lambda}) + \lambda^2 e^{-\lambda} (2-e^{-\lambda}) \nonumber
\end{align}
where $\tilde{\nu}_i^\prime, i\in\Sigma$, denotes the derivative of 
$\nu_i$, i.e. for  $z>0$
\begin{align*}
\tilde{\nu}_i^\prime (z)=\sum_{n=1}^\infty e^{-z} 
\frac {z^{n-1}}{(n-1)!} \E[L_n^i]-\tilde{\nu}_i (z).
\end{align*}
\end{lem}

The Mellin transform $v_i^*(s)$ of $\tilde{v}_i (\lambda)$
is
\begin{align*}
v_i^*(s)&=p_{i0}^{-s}v_0^*(s)+p_{i1}^{-s}v_1^*(s)
-2s p_{i0}^{-s} \nu_0^*(s) \\
&\quad -2s p_{i1}^{-s} \nu_1^*(s) -\Gamma(s+1)+F_i^*(s)
\end{align*}
with $F_i^*(s)$  the Mellin transform of 
$e^{-\lambda}(
\tilde{\nu}_0^\prime 
(\lambda p_{i0}) + 2 \lambda p_{i1} \tilde{\nu}_1^\prime (\lambda p_{i1})
+\lambda^2 (2-e^{-\lambda}))$. Thus, the column vector $\bv^*(s):=(v_0^*(s), 
v_1^*(s))$ satisfies the following algebraic equation
\begin{align*}
\bv^*(s)=& P(s)\bv^*(s)-2 s P(s-1)\vnu^*(s)\\
&~-\vgamma(s+1)+{\bf F^*}(s)
\end{align*}
where ${\bf F^*}(s):=(F^*_0(s), F^*_1(s))$. 
Then, as we did before for the mean analysis, we obtain
\begin{align*}
\bv(s)&=- \frac{2s \Gamma(s+1) \la \vpi(s), P(s-1) \VPsi(s) \ra 
\VPsi(s)}{(1-\lambda(s))^2} \\
&\quad  + \bo(1/(1-\lambda(s)).
\end{align*}
After further  computations we find that the Poisson variance 
$\tilde{v}(\lambda)={\rm Var}(L^\mu_{N_\lambda})$ is
\begin{align*}
\tilde{v}(\lambda)&=\frac{1}{\dlam^2(-1)} \lambda  \log^2 \lambda +
\left(\frac{\ddlam(-1)}{2\dlam^3(-1)}+\frac{A}{\dlam^2(-1)} \right)
\lambda \log \lambda \\
&\quad +\bo(\lambda)
\end{align*}
for some explicitly computable constant $A$. Finally, with
depoissonization, cf.~\cite{spa-book}, we obtain
\begin{align*}
{\rm Var}(L^\mu_n)&= \tilde{v}(n)-n [\tilde{\nu}^\prime(n)]^2 \\
&=\frac{\ddlam(-1)-\dlam^2(-1)}{\dlam^3(-1)} n \log n +\bo(n)
\end{align*}
proving Theorem~\ref{thm:variance}.

\section{Asymptotic Normality} \label{sec:limit}

Our main result is the asymptotic normality of the external path length:

\begin{thm} \label{thm:limit}
For the external path length $L_n^\mu$ of a binary trie under the Markov source model with conditions (\ref{cond_prob}) we have 
\begin{align}\label{limit_law}
 \frac{L_n^\mu - \E[L_n^\mu]}{\sqrt{n\log n}} \stackrel{d}{\longrightarrow} {\cal N}(0,\sigma^2), \qquad (n\to\infty),
\end{align}
where $\sigma^2>0$ is independent of the initial distribution $\mu$ and given by (\ref{var_formular}).
\end{thm}

As in the analysis of the mean, we first derive limit laws for $L_n^0$ and $L_n^1$ and then transfer these to a limit law for $L_n^\mu$ via (\ref{decomp_initial}). We abbreviate for $i\in\Sigma$ and $n\in\N_0$
\begin{align*}
\nu_i(n):=\E[L_n^i],   \qquad \sigma_i(n):=\sqrt{\V(L_n^i)}.
\end{align*}
Note that we have $\nu_i(0)=\nu_i(1)=\sigma_i(0)=\sigma_i(1)=0$ and $\sigma_i(n)>0$ for all $n\ge 2$. 
We define the standardized variables by 
\begin{align}\label{normal_rec}
Y^i_n := \frac{L_n^i - \E[L_n^i]}{\sigma_i(n)}, \qquad i\in\Sigma, n\ge 2,
\end{align}
and $Y^i_0:=Y^i_1:=0$. Then we have:

\begin{prop} \label{conv_zolo_prop}
For both sequences $(Y^i_n)_{n\ge 0}$, $i\in\Sigma$, we have convergence 
in distribution:
\begin{align} \label{conv_zolo}
Y^i_n \stackrel{d}{\longrightarrow}{\cal N}(0,1) \qquad (n\to\infty).
\end{align}
\end{prop}

We now present a brief streamlined road map of the proof.
\med
{\bf Step 1. Normalization.} From the system (\ref{decomp_zero})--(\ref{decomp_one}), where we denote there   $I^0_n:=I_n$ and $I^1_n:=J_n$, 
and the normalization (\ref{normal_rec})  we obtain for all $n\ge 2$,
\begin{align}\label{mod_rec}
Y^i_n \stackrel{d}{=}\frac{\sigma_i(I_n^i)}{\sigma_i(n)} Y^i_{I_n^i}+ \frac{\sigma_{1-i}(n-I_n^i)}{\sigma_i(n)}Y^{1-i}_{n-I_n^i} + b_i(n), 
\end{align}
where 
\begin{align*}
 b_i(n)= \frac{1}{\sigma_i(n)}\left( n +\nu_i(I_n^i) + 
\nu_{1-i}(n-I^i_n)-\nu_i(n) \right),
\end{align*}
and in (\ref{mod_rec}) we have that $(Y^0_0,\ldots,Y^0_n)$, 
$(Y^1_0,\ldots,Y^1_n)$ and 
$(I_n^0,I_n^1)$ are independent.  It can be shown by our expansions of the 
means $\nu_i(n)$ and the Lipschitz property from Proposition \ref{thEwLip} 
that we have $b_i(n) \to 0$ as $n\to\infty$ for both $i\in \Sigma$, e.g., in 
the $L_3$-norm which below will be technically sufficient.  
Furthermore, the asymptotic of the variance from 
Theorem~\ref{thm:variance} implies together with the strong law of 
large numbers that the coefficients in (\ref{mod_rec}) converge:
 \begin{align*}
\frac{\sigma_i(I_n^i)}{\sigma_i(n)}  \to \sqrt{p_{ii}}, \qquad   
\frac{\sigma_{1-i}(n-I_n^i)}{\sigma_i(n)}\to \sqrt{1-p_{ii}},
\end{align*}
where we recall that $\sigma_i(I_n^i)$ is the standard deviation of $L^i_{I_n^i}$ conditioned on
$I_n^i$, hence, in particular a random variable. \\

\noindent
{\bf Step 2. System of limit equations.}
The convergence of the coefficients in (\ref{mod_rec}) suggests, by passing formally with $n\to\infty$,  that limits $Y^0$ and $Y^1$ of  
$Y^0_n$ and $Y^1_n$, if they exist, should satisfy the 
system of recursive distributional equations
\begin{align}\label{sys_fix1}
Y^0 &\stackrel{d}{=} \sqrt{p_{00}} Y^0 + \sqrt{1-p_{00}} Y^1 ,\\
Y^1 &\stackrel{d}{=} \sqrt{1-p_{11}} Y^0 + \sqrt{p_{11}} Y^1, \label{sys_fix2}
\end{align}
where $Y^0$ and $Y^1$ are being independent on the right hand sides. 
Clearly, centered normally distributed $Y^0$ and $Y^1$ with 
identical variances solve the system (\ref{sys_fix1})--(\ref{sys_fix2}). The task now is to show that $Y^0_n$ and $Y^1_n$  converge in distribution towards these solutions $Y^0$ and $Y^1$ respectively.\\

\noindent
{\bf Step 3. The operator of  distributions.}
Our approach is based on the system (\ref{sys_fix1})--(\ref{sys_fix2}) of limit equations  together with an associated 
contracting operator (map) on the space of probability distributions as follows: 
We denote by ${\cal M}_s(0,1)$ the space of all probability distributions 
on the real line with mean $0$, variance $1$ and finite absolute  
moment of order $s$. Later $2<s\le 3$ will be an appropriate 
choice for us. With the abbreviation 
${\cal M}^2:= {\cal M}_s(0,1) \times {\cal M}_s(0,1)$  we define the map
\begin{align*}
T:  {\cal M}^2 &\to {\cal M}^2\\
(\tau_0,\tau_1) &\mapsto \left({\cal L}\left(\sqrt{p_{00}} 
W^0 + \sqrt{1-p_{00}} W^1\right)\right., \\
&\qquad \left. {\cal L}\left(\sqrt{1-p_{11}} W^0 + \sqrt{p_{11}} W^1\right)\right),
\end{align*}
where $W^0$, $W^1$ are independent with distributions 
${\cal L}(W^i)=\tau_i$ for both $i\in\Sigma$. 

This allows a measure theoretic reformulation of solutions of (\ref{sys_fix1})--(\ref{sys_fix2}) that is convenient subsequently:
Random variables $(Y^0,Y^1)$ solve the system (\ref{sys_fix1})--(\ref{sys_fix2}) if and only if their 
pair of distributions  $({\cal L}(Y^0),{\cal L}(Y^1))$ is a fixed point of $T$. Hence the identification of fixed-points and domains of attraction of such fixed-points plays an important role in the asymptotic behavior of our sequences $(Y_n^0)_{n\ge 0}$ and $(Y_n^1)_{n\ge 0}$ and is a core part of our proof.\\

\noindent
{\bf Step 4. The Zolotarev metric.} In accordance with the general idea of the contraction method we will endow the space ${\cal M}^2$ with a complete metric such that $T$ becomes a contraction with respect to this metric. The issue of fixed-points is then reduced to the application of   Banach's fixed-point theorem.

As building block we use the 
{\it Zolotarev metric} on  ${\cal M}_s(0,1)$.  It  has been 
studied in the context of the contraction method systematically 
in \cite{NeRu04}. We only need
the following properties, see Zolotarev \cite{zo76,zo77}:
For distributions 
${\cal L}(X)$, ${\cal L}(Y)$ on $\R$ the Zolotarev 
distance $\zeta_s$, $s>0$, is defined by
\begin{align}
\label{eq:3.6} \zeta_s(X,Y) &:= \zeta_s({\cal L}(X),{\cal L}(Y))\\
&:=\sup_{f\in {\cal F}_s}|\E[f(X) -
f(Y)]|\nonumber
\end{align}
where $s=m+\alpha$ with $0<\alpha\le 1$,
$m\in\N_0$, and
\begin{align*} 
{\cal F}_s:=\{f\in
C^m:\|f^{(m)}(x)-f^{(m)}(y)\|\le
\|x-y\|^\alpha\},
\end{align*}
 the space of $m$ times
continuously differentiable functions from
$\R$ to $\R$ such that the $m$-th
derivative is H\"older continuous of order
$\alpha$ with H\"older-constant $1$. 
We have that $\zeta_s(X,Y)<\infty$, if all 
moments of orders $1,\ldots,m$ of $X$ and $Y$ are
equal and if the $s$-th absolute moments of $X$ and
$Y$ are finite.  Since later on only the case $2<s\le 3$ is
used, for finiteness of  $\zeta_s(X,Y)$ it is thus sufficient  for these 
$s$ that  mean and variance of 
$X$ and $Y$ coincide and both have a finite absolute moment of order $s$.
Convergence in $\zeta_s$ implies weak convergence on $\R$.
Furthermore,  $\zeta_s$ is $(s,+)$ ideal, i.e., we have
\begin{align*}
&\zeta_s(X+Z,Y+Z)\le\zeta_s(X,Y), \\
 &\zeta_s(cX,cY) = c^s \zeta_s(X,Y)
\end{align*}
for all  $Z$ being independent of $(X,Y)$ and all $c>0$.

Now, to measure distances on the product space ${\cal M}^2$ we define for $(\tau_0,\tau_1), (\varrho_0,\varrho_1) \in {\cal M}^2$ the distance
\begin{align*}
\zeta_s^\vee((\tau_0,\tau_1), (\varrho_0,\varrho_1)):=\zeta_s(\tau_0,\varrho_0) \vee \zeta_s(\tau_1,\varrho_1).
\end{align*}
Here and later on, we use the symbols $\vee$ and $\wedge$ for $\max$ and $\min$ respectivly.

\noindent
{\bf Step 5. The contraction property}.
We directly obtain that $T$ is a contraction in $\zeta_s^\vee$ from the property that $\zeta_s$ is $(s,+)$ ideal: 
Denoting the components of $T$ by $T_0$ and $T_1$ we have 
\begin{align*}
&\zeta_s(T_0(\tau_0,\tau_1), T_0(\varrho_0,\varrho_1)) \\
&\le p_{00}^{s/2} \zeta_s(\tau_0,\varrho_0) + (1-p_{00})^{s/2} \zeta_s(\tau_1,\varrho_1) \\
& \le \left( p_{00}^{s/2} + (1-p_{00})^{s/2}\right) \zeta_s^\vee((\tau_0,\tau_1), (\varrho_0,\varrho_1)),
\end{align*}
and similary
\begin{align*}
&\zeta_s(T_1(\tau_0,\tau_1), T_1(\varrho_0,\varrho_1))\\
&\le (1-p_{11})^{s/2} \zeta_s(\tau_0,\varrho_0) + p_{11}^{s/2} \zeta_s(\tau_1,\varrho_1) \\
&\le \left( (1-p_{11})^{s/2} + p_{11}^{s/2}\right) \zeta_s^\vee((\tau_0,\tau_1), (\varrho_0,\varrho_1)).
\end{align*}
Hence together with $\xi:= \max_{i\in\Sigma} (p_{ii}^{s/2} + (1-p_{ii})^{s/2})$ we obtain that
\begin{align}\label{nnn_cont}
\zeta_s^\vee(T(\tau_0,\tau_1), T(\varrho_0,\varrho_1))\le \xi  \zeta_s^\vee((\tau_0,\tau_1), (\varrho_0,\varrho_1)).
\end{align}
Since $p_{ii}\in (0,1)$ by assumption (\ref{cond_prob}) we have $\xi<1$ for all $s>2$.  On the other hand, it is known that one only obtains finiteness of
$\zeta_s$ on ${\cal M}_s(0,1)$ for $s\le 3$, hence (\ref{nnn_cont}) is only meaningful for  $s\le 3$. Thus, altogether, our choice of $s$ is 
$2<s\le 3$. For these $s$ we obtain that $T$ is a contraction in 
$\zeta_s^\vee$.  \\

\noindent
{\bf Step 6. Convergence of the }$\mathbf{Y^i_n}$.
An intuition why contraction properties of   the map $T$ lead to  
convergence of the $Y^i_n$ towards the unique fixed-point $({\cal N}(0,1), {\cal N}(0,1))$ of $T$ in ${\cal M}^2$ is as follows:  The map $T$ serves as a limit version of our 
recurrence system (\ref{mod_rec}). Since in this recurrence system we could  replace the $Y^i_{I^i_n}$ and  $Y^{1-i}_{n-I^i_n}$ on the right hand side by the recurrence (\ref{mod_rec}) itself, iterating these replacements  leads approximatively to an iteration of the map $T$. However, by Banach's fixed-point theorem, the iteration of $T$ applied to any starting point in  ${\cal M}^2$ converges to the unique fixed-point of $T$ in the metric $\zeta_s^\vee$. 

Hence, the problem of proving the convergence of the $Y^i_n$ to the standard normal distribution (the fixed-point) is reduced to the following technical task: Verify that not only the iterations of $T$ itself convergence in the metric $\zeta_s^\vee$ to the fixed-point, but also that the iterations of the approximations of $T$ that make the recurrence of the  $Y^i_n$ convergence within  $\zeta_s^\vee$. 

Once this is settled, we use that convergence in $\zeta_s$ is strong enough to imply weak convergence and $({\cal N}(0,1), {\cal N}(0,1))$ is the unique fixed point of $T$. This finally yields Proposition \ref{conv_zolo_prop}.
A detailed proof is given in the full paper version of this extended abstract.\\
\med
{\bf Step 7. Transfer to arbitrary initial distributions.} 
Finally, we prove Theorem \ref{thm:limit}. For this, we have to transfer the convergence of the 
$Y^i_n$ from Proposition \ref{conv_zolo_prop} to the convergence of the normalization of $L^\mu_n$ via (\ref{decomp_initial}). Recall that in (\ref{decomp_initial}), the $K_n$ is a binomial $B(n,\mu_0)$ distributed random variable. 
We write 
\begin{align*}
\frac{L_n^\mu - \E[L_n^\mu]}{\sqrt{n\log n}}
&=\frac{L_n^\mu - \nu_0(K_n)- \nu_1(n-K_n)}{\sqrt{n\log n}} \\
&\quad + \frac{\nu_0(K_n)+ \nu_1(n-K_n) - \E[L_n^\mu]}{\sqrt{n\log n}}.
\end{align*}
By the Lemma of Slutsky, see, e.g.~\cite[Theorem 3.1]{bi99}, it is sufficient to show, as $n\to\infty$,
\begin{align}
\frac{L_n^\mu - \nu_0(K_n)- \nu_1(n-K_n)}{\sqrt{n\log n}} & \stackrel{d}{\longrightarrow} {\cal N}(0,\sigma^2)  \label{conv_1}\\
 \frac{\nu_0(K_n)+ \nu_1(n-K_n) - \E[L_n^\mu]}{\sqrt{n\log n}} & \stackrel{\Prob}{\longrightarrow} 0.\label{conv_2}
\end{align}
For showing (\ref{conv_1}) note that by Proposition \ref{conv_zolo_prop} $(L_n^i - \E[L_n^i])/\sqrt{n\log n} \to {\cal N}(0,\sigma^2) $ in distribution for both $i\in\Sigma$. 
We set $A_n:= [\mu_0 n -n^{2/3}, \mu_0 n +n^{2/3}] \cap \N_0$ and $A_n^c:=\{0,\ldots,n\}\setminus A_n$. Then by Chernoff's bound (or the central limit theorem) we have $\Prob(K_n \in A_n) \to 1$. For all $x\in\R$ we have with $\kappa_{nj}:=\Prob (K_n=j)$
\begin{align*}
\lefteqn{\Prob\left(\frac{L_n^\mu - \nu_0(K_n)- \nu_1(n-K_n)}{\sqrt{n\log n}} \le x\right)}\\
&=\Prob\left(\frac{L_{K_n}^0 - \nu_0(K_n)}{\sqrt{n\log n}} +  \frac{L_{n-K_n}^1   - \nu_1(n-K_n)}{\sqrt{n\log n}} \le x\right)  \\
&= \sum_{j\in A_n}\kappa_{nj}\Prob\left(\frac{L_{j}^0 - \nu_0(j)}{\sqrt{n\log n}} +  \frac{L_{n-j}^1   - \nu_1(n-j)}{\sqrt{n\log n}} \le x\right) \\
&\quad +o(1).
\end{align*}
For $j\in A_n$ we have $\sqrt{j\log j}/\sqrt{n\log n} \to \sqrt{\mu_0}$ and 
$\sqrt{(n-j)\log(n- j)}/\sqrt{n\log n} \to \sqrt{1-\mu_0}$. Hence, we have $(L_{j}^0 - \nu_0(j))/\sqrt{n\log n} \to 
{\cal N}(0,\mu_0 \sigma^2)$ and $(L_{n-j}^1   - \nu_1(n-j))/\sqrt{n\log n}\to 
{\cal N}(0,(1-\mu_0) \sigma^2)$ in distribution and the two summands are independent. Together, denoting by $N_{0,\sigma^2}$ an  ${\cal N}(0,\sigma^2)$ distributed random variable we obtain
\begin{align*}
&\Prob\left(\frac{L_n^\mu - \nu_0(K_n)- \nu_1(n-K_n)}{\sqrt{n\log n}} \le x\right) \\
&= o(1) + \sum_{j\in A_n} \kappa_{nj} ( \Prob\left(N_{0,\sigma^2} \le x\right) +o(1))\\
&\to \Prob\left(N_{0,\sigma^2} \le x\right),
\end{align*}
where the latter convergence is justified by dominated convergence. This shows (\ref{conv_1}).

To establish the  convergence in probability in (\ref{conv_2}) note that 
(\ref{decomp_initial}) implies
$$\E[L_n^\mu] = \E[\nu_0(K_n)] + \E[\nu_1(n-K_n)].$$
Hence, with the notation (\ref{def_fi}) 
and $g(x):=x\log x$ for $x\in[0,1]$ and $\| \, \cdot\,\|_1$ denoting the $L_1$-norm we have
\begin{align*}
\lefteqn{\frac 1 {\sqrt{n\log n}} \| \nu_0(K_n)+\nu_1(n-K_n)- \E[L_n^\mu] \|_1}  \\
&=\frac 1 {\sqrt{n\log n}} \| \nu_0(K_n)-\E[\nu_0(K_n)] \\
& \hspace{1.9cm}+\nu_1(n-K_n)-\E[\nu_1(n-K_n)] \|_1 \\
&\leq \frac 1 {H \sqrt{n\log n}} \| g(K_n)-\E[g(K_n)]\\
&\hspace{2.2cm} +g(n-K_n)-\E[g(n-K_n)] \|_1 \\
&\;\;\;\;~+ \frac 1 {\sqrt{n\log n}} \| f_0(K_n)-\E[f_0(K_n)] \|_1 \\
&\;\;\;\;~+\frac 1 {\sqrt{n\log n}} \| f_1(n-K_n)-\E[f_1(n-K_n)] \|_1.
\end{align*}
With the concentration of the binomial distribution we obtain
\begin{align*}
\lefteqn{\| g(K_n)-\E[g(K_n)]+g(n-K_n)-\E[g(n-K_n)] \|_1 } \quad\phantom{.}\\
&=n\left\|  g\left( \frac{K_n}{n} \right)-\E \left[g\left(\frac{K_n}{n}\right)\right] \right.\\
&\qquad \left.+g\left(\frac{n-K_n}{n}\right)-\E\left[ g\left(\frac{n-K_n}{n} \right)\right] \right\|_1\\
&= \bo\left(n^{\nicefrac{1}{2}}\right).
\end{align*}
The terms $\| f_0(K_n)-\E[f_0(K_n)] \|_1$ and $\| f_1(n-K_n)-\E[f_1(n-K_n)] \|_1$ are also of the order $\bo(n^{\nicefrac{1}{2}})$ by a self-centering argument.
 Altogether we have
$$\frac{\| \nu_0(K_n)+\nu_1(n-K_n)- \E[L_n^\mu] \|_1}{\sqrt{n\log n}}= \bo\left( \frac{1}{\sqrt{\log n}}\right), $$
which, by Markov's inequality, implies (\ref{conv_2}) as follows: For any $\varepsilon > 0$ we have
\begin{align*}
\lefteqn{\Prob\left( \left|\frac{\nu_0(K_n)+ \nu_1(n-K_n) - \E[L_n^\mu]}{\sqrt{n\log n}}\right| > \varepsilon\right) }\\
&\le \frac{1}{\varepsilon} \E \left[
 \left|\frac{\nu_0(K_n)+ \nu_1(n-K_n) - \E[L_n^\mu]}{\sqrt{n\log n}}\right|\right] \\
&=\frac{1}{\varepsilon\sqrt{n\log n}}  \| \nu_0(K_n)+\nu_1(n-K_n)- \E[L_n^\mu] \|_1 \\
&\to 0.
\end{align*}

\section{Comparison with a multivariate approach}
\label{rem_multivariate}
We propose the use of systems of univariate recurrences in this extended abstract. Note however, that known limit theorems  from the contraction 
method for multivariate recurrences can as well be applied to the bivariate random variable 
$Y_n:=(Y^0_n,Y^1_n)$. (Technically easiest is to keep the components $Y^0_n$ and $Y^1_n$  independent by working with independent $I_n^0$ and $I_n^1$.)
Applying such an approach as developed in \cite{NeRu04},  the system (\ref{sys_fix1})--(\ref{sys_fix2}) is now replaced by the bivariate recursive distributional equation 
\begin{align} \label{limit_2d_app}
Y \stackrel{d}{=}A_1Y+ A_2 \widehat{Y}, 
\end{align}
where $Y$ and $\widehat{Y}$ are independent and identically distributed bivariate random variables and the matrices $A_1, A_2$ are give by 
\begin{align*}
A_1 &:=  \left[ \begin{array}{cc} \sqrt{p_{00}} & 0 \\  0 & \sqrt{p_{11}} \end{array} \right],\\
A_2 &:= \left[ \begin{array}{cc} 0 & \sqrt{1-p_{00}}  \\   \sqrt{1-p_{11}} & 0\end{array} \right].
\end{align*}
Any centered bivariate normal distribution solves the latter fixed-point equation (\ref{limit_2d_app}). In particular Theorem 4.1 in \cite{NeRu04} covers the arising bivariate recurrence, cf.~also condition (38) in \cite{NeRu04}, which is satisfied for $A_1$, $A_2$  in (\ref{limit_2d_app})

However, for applying the contraction method in such a multivariate form, an underlying contraction
 is only implied for, see condition (25) in \cite{NeRu04},
\begin{align*}
 \|A_1\|_\mathrm{op}^3 + \|A_2\|_\mathrm{op}^3<1,
\end{align*}
where $\|\cdot\|_\mathrm{op}$, here, is identical to the spectral radius of the matrix.
This imposes the additional condition 
\begin{align} \label{cond_prob_add_app}
 (p_{00}\vee p_{11})^{3/2} + (1 - p_{00}\wedge p_{11})^{3/2} < 1
\end{align}
to come up with a result similar to our Theorem \ref{thm:limit}.

Our new approach based on systems of  univariate recursive equations given above does not require any further condition such as (\ref{cond_prob_add_app}).

\end{document}